\documentclass{amsart}
\usepackage[english]{babel}
\usepackage[latin1]{inputenc}
\usepackage[dvips,final]{graphics}
\usepackage{amsmath,amsfonts,amssymb,amsthm,amscd,array,stmaryrd,mathrsfs}
\usepackage{pstricks}
 \usepackage[all]{xy}
\usepackage{textcomp}
 \usepackage[final]{epsfig}
\theoremstyle{plain}
\newtheorem{thm}{Theorem}%[section]
\newtheorem{lem}{Lemma}[section]
\newtheorem{cor}[lem]{Corollary}
\newtheorem{prop}[lem]{Proposition}

\theoremstyle{definition}
\newtheorem{defn}[lem]{Definition}
\newtheorem{rem}[lem]{Remark}
\newtheorem{ex}[lem]{Example}

\newcommand{\R}{\mathbb{R}}
\newcommand{\Z}{\mathbb{Z}}
\newcommand{\C}{\mathbb{C}}
\newcommand{\N}{\mathbb{N}}

\newcommand{\K}{\mathbb{K}}

\newcommand{\Dc}{\mathcal{D}}

\newcommand{\F}{\mathcal{F}}
\newcommand{\cc}{\Gamma}

\newcommand{\E}{\mathcal{E}}

\newcommand{\gog}{\mathfrak{g}}

\newcommand{\ga}{\mathfrak{a}}
\newcommand{\U}{\mathcal{U}}

\newcommand{\id}{\textup{Id}}

\newcommand{\half}{\frac{1}{2}}
\newcommand{\End}{\textup{End}}
\newcommand{\Hom}{\textup{Hom}}
\newcommand{\e}{\varepsilon}
\newcommand{\ad}{\mathrm{ad}}

\def\a{\alpha}

\def\e{\varepsilon}

\def\l{\lambda}

\begin{document}

\title{Representations of $\mathrm{asl}_2$}

\author{Sophie Morier-Genoud}

\address{Department of Mathematics, University of Michigan,
Ann Arbor, USA}

\email{sophiemg@umich.edu}

\date{}


\begin{abstract}
We study representations of the simple Lie
antialgebra $\mathrm{asl}_2$ introduced in \cite{Ovs}.
We show that representations of $\mathrm{asl}_2$ are
closely related to the famous ghost Casimir element
of the universal enveloping algebra $\mathrm{osp}(1|2)$.
We prove that $\mathrm{asl}_2$ has no non-trivial finite-dimensional
representations; we define and classify some particular
infinite-dimensional representations that we call weighted representations.
\end{abstract}

\maketitle

\vspace{-0.8cm}
%\tableofcontents

%%%%%%%%%%%%%%%%%%%%%%
\section*{Introduction}
%%%%%%%%%%%%%%%%%%%%%%

``Lie antialgebras'' is a new class of algebras introduced
by V. Ovsienko \cite{Ovs}.
These algebras naturally appear in the context of symplectic and contact
geometry of $\Z_2$-graded spaces, their algebraic properties
are not yet well understood.
Lie antialgebras is a surprizing ``mixture'' of commutative algebras and
Lie algebras.

The first example of Ovsienko's algebras is a simple Lie antialgebra
$\mathrm{asl}_2(\K)$, over $\K=\R$ or $\C$.
This algebra is of dimension 3, it has linear basis $\{\e;\,a,b\}$
subject to the following relations:
\begin{equation}
\label{aslA}
\begin{array}{llllllllll}
\e\cdot\e&=&\e,&&&&&&&
\\[4pt]
\e\cdot{}a&=&a\cdot{}\e&=&\half \,a,&
\qquad
\e\cdot{}b&=&b\cdot{}\e&=&\half \,b,
\\[5pt]
a\cdot{}b&=&-b\cdot{}a&=&\half \,\e.&&&&&
\end{array}
\end{equation}
The notion of representation of a Lie antialgebra was also defined
in \cite{Ovs}, and the problem of classification of representations
of simple Lie antialgebras was formulated.
In this paper, we study representations of $\mathrm{asl}_2(\K)$.

The Lie antialgebra $\mathrm{asl}_2(\K)$ is closely related to
the simple classical Lie superalgebra $\mathrm{osp}(1|2)$.
For instance, 
$\mathrm{osp}(1|2)=\mathrm{Der}\left(\mathrm{asl}_2(\K)\right)$.
It was shown in \cite{Ovs}, that every representation of
$\mathrm{asl}_2(\K)$ is naturally a representation of $\mathrm{osp}(1|2)$.
The first problem is thus to determine the corresponding class of
$\mathrm{osp}(1|2)$-representations.
It turns out that this class can be characterized by so-called ghost Casimir element
of the universal enveloping algebra $\U(\mathrm{osp}(1|2))$. This element first appeared in \cite{Pinczon}, see also
\cite{ABP, Fra, Gor}.

\begin{thm}
\label{NoOne}
There is a one-to-one correspondence between representations
of $\mathrm{asl}_2(\K)$ and representations of
$\mathrm{osp}(1|2)$ such that 
\begin{equation}
\label{GhsAct}
\textstyle
\cc^2=\frac{1}{4}\,\id,
\end{equation}
where
$\cc$ is the action of the ghost Casimir element.
\end{thm}

We will see that in the case of a $\Z_2$-graded irreducible representation $V=V_0\oplus{}V_1$
of $\mathrm{asl}_2(\K)$, one has a more complete information on the
action of the ghost Casimir element, namely, the $\Z_2$-grading can be chosen
in such a way that
$$
\textstyle
\cc\left|_{V_0}\right.=-\frac{1}{2}\,\id,
\qquad
\cc\left|_{V_1}\right.=\frac{1}{2}\,\id.
$$
We believe that this relation with the ghost Casimir element
provides a better understanding for the nature of $\mathrm{asl}_2(\K)$ itself.

In the finite-dimensional case, we prove the following
\begin{thm}
\label{FinThm}
The Lie antialgebra $\mathrm{asl}_2(\K)$ has no
non-trivial finite-dimensional representations.
\end{thm}
\noindent
Let us emphasize that
the algebra $\mathrm{asl}_2(\K)$ was defined in \cite{Ovs}
as analog of the classical simple Lie algebra $\mathrm{sl}_2(\K)$,
but is also has a certain similarity with the 3-dimensional
Heisenberg algebra $\mathfrak{h}_1$.
The above result is similar to the classical result
that $\mathfrak{h}_1$ has no non-trivial irreducible
representations.

In the infinite-dimensional case, we restrict our study to the case of  \textit{weighted representations}.
A weighted representation is a representation 
containing an eigenvector for the action of the Cartan element $H\in\mathrm{osp}(1|2)$. 
We classify the irreducible weighted representations of $\mathrm{asl}_2(\K)$. We introduce a family of weighted representations $V(\ell)$, for $\ell \in \K$ (see Section \ref{Vell} for the construction). Considering the set of parameters $\mathcal{P}=[-1,1]$ in the real case, or $\mathcal{P}=[-1,1]\cup \{\ell \in \C\; |\; -1\leq \mathrm{Re}(\ell) <1\}$ in the complex case, we obtain the complete classification of irreducible weighted representations.

\begin{thm}\label{InfThm}
Any irreducible weighted representation is isomorphic to a  $V(\ell)$ for a unique 
$\ell \in \mathcal{P}$.
\end{thm}

The paper is organized as follows. In Section 1, we recall the general definitions of Lie antialgebras and their representations introduced in \cite{Ovs}. In Section 2, we obtain preliminary results about the representations of $\mathrm{asl}_2$. A link with the Lie algebra $\mathrm{osp}(1|2)$ and with the Casimir elements is established. We complete the proof of Theorem \ref{NoOne} in subsection \ref{PfNoOne}. In Section 3, we introduce the notion of weighted representations and give the construction of the family of irreducible weighted representations $V(\ell)$, $\ell \in \K$. In Section 4, we formulate our results concerning the representations $V(\ell)$ and complete the proofs of Theorem \ref{FinThm} and Theorem \ref{InfThm}. In the end of the paper we discuss some general aspects of
representation theory of $\mathrm{asl}_2(\K)$, such
as the tensor product of two representations. 

\bigskip
\noindent
\textbf{Acknowledgements}.
I am grateful to V. Ovsienko for the statement of the
problem and enlightening discussions.

%%%%%%%%%%%%%%%%%%%%%%
\section{Lie antialgebras and their representations}
%%%%%%%%%%%%%%%%%%%%%%

Let us give the definition of a Lie antialgebra
equivalent to the original definition of \cite{Ovs}. Throughout the paper the ground vector field is $\K=\R$ or $\C$.

\begin{defn}
A Lie antialgebra is
a $\Z_2$-graded vector space $\ga=\ga_0\oplus\ga_1$,  
equipped with a bilinear product $\cdot$ satisfying the following conditions.

\begin{enumerate}
\item 
it is even:
$
\ga_i\cdot\ga_j\subset\ga_{i+j};
$

\item 
it is supercommutative, \textit{i.e.},
for  all homogeneous elements $x, y \in \ga$,
\begin{equation*}
\label{SkewP}
x\cdot{}y=\,(-1)^{p(x)p(y)}\,
y\cdot{}x
\end{equation*}
where $p$ is the parity function defined by $p(x)=i$
for $x\in \ga_i$;

\item
the subspace $\ga_0$ is a
commutative associative algebra;

\item
for all $x_1,x_2\in\ga_0$ and $y\in\ga_1$, one has
$$
\textstyle
x_1\cdot\left(x_2\cdot{}y\right)=
\frac{1}{2}\left(x_1\cdot{}x_2\right)\cdot{}y,
$$
in other words,
the subspace $\ga_1$ is a module over $\ga_0$,
homomorphism
 $\varrho:\ga_0\to\End(\ga_1)$ being given by
 $
 \varrho_x\,y=2\,x\cdot{}y
 $
 for all $x\in\ga_0$ and $y\in\ga_1$;

\item
for all $x\in\ga_0$ and $y_1,y_2\in\ga_1$, the following Leibniz
identity
$$
\textstyle
x\cdot\left(y_1\cdot{}y_2\right)=
\left(x\cdot{}y_1\right)\cdot{}y_2+
y_1\cdot\left(x\cdot{}y_2\right)
$$
is satisfied;

\item
for all $y_1,y_2,y_3\in\ga_1$, the following Jacobi-type
identity
$$
y_1\cdot\left(y_2\cdot{}y_3\right)+
y_2\cdot\left(y_3\cdot{}y_1\right)+
y_3\cdot\left(y_1\cdot{}y_2\right)=0
$$
is satisfied.
\end{enumerate}
\end{defn}

\begin{ex}
It is easy to see that the above axioms are satisfied
for $\mathrm{asl}_2(\K)$.
In this case,
the element $\e$ spans the even part, $\mathrm{asl}_2(\K)_0$,
while the elements $a,b$ span the odd part,
$\mathrm{asl}_2(\K)_1$.
\end{ex}

Consider a $\Z_2$-graded vector space  $V=V_0\oplus{}V_1$,
the space $\End(V)$ of linear endomorphisms of $V$ is
a $\Z_2$-graded associative algebra:
$$
\End(V)_0=
\End(V_0)\oplus\End(V_1),
\qquad
\End(V)_1=
\Hom(V_0,V_1)\oplus
\Hom(V_1,V_0).
$$
Following \cite{Ovs}, we define the following
``anticommutator'' on $\End(V)$:
\begin{equation}
\label{AntiCoCo}
\left]X,Y\right[:=
X\,Y+(-1)^{p(X)p(Y)}
Y\,X,
\end{equation}
where $p$ is the parity function on $\End(V)$
and  $X,Y\in\End(V)$ are homogeneous (purely even or purely odd)
elements.
Note that the sign rule in (\ref{AntiCoCo}) is opposite
to that of the usual commutator.

\begin{rem}
Let us stress that the operation (\ref{AntiCoCo}) 
does \textit{not} define a Lie antialgebra structure
on the full space $\End(V)$ and it is not known what are the
subspaces of $\End(V)$ for which this is the case.
This operation provides, however, a definition
of the notion of representation of a Lie antialgebra.
\end{rem}

\begin{defn}
\label{defRep}
(a)
We call \textit{representation} of the Lie antialgebra  $\ga$,  the data of $(V, \chi)$ where $V=V_0\oplus{}V_1$ is  a $\Z_2$-graded vector space
and $\chi:\ga\to\End(V)$ is an even linear map
such that
\begin{equation}\label{eqRep}
\left]\chi_x,\chi_y\right[=
\chi_{x\cdot{}y},
\end{equation}
for all $x,y\in\ga$.

(b)
A \textit{subrepresentation} is a $\Z_2$-graded subspace $V' \subset V$ stable
under $\chi_x$ for all $x\in\ga$.

(c)
A representation is called \textit{irreducible} if it does not have 
proper subrepresentations.

(d)
Two representations $(V,\chi)$ and $(V',\chi')$ are called \textit{equivalent} if there exists a  linear map
$\Phi:V\to{}V'$ such that $\Phi \circ \chi_x=\chi'_x \circ \Phi $, for every $x\in \ga$.
\end{defn}

\begin{rem}
It is not clear \textit{a priori}, that a given Lie antialgebra has at least
one non-trivial representation.
In the case of $\mathrm{asl}_2(\K)$, however, an example
of representation was given in \cite{Ovs} in the context of
geometry of the supercircle.

Let us finally mention that there is a notion of module
over a Lie antialgebra which is different from that of
representation.
For instance, the ``adjoint action'' defined as usual by
$\ad_x\,y=x\cdot{}y$ is \textit{not} a representation,
but it defines a structure of $\ga$-module on $\ga$.
\end{rem}

Given a Lie antialgebra $\ga$, it was shown in \cite{Ovs} that there exists
a Lie superalgebra, $\gog_\ga$, canonically associated to $\ga$.
Every representation of $\ga$ extends to a representation of
$\gog_\ga$.
In the case of $\mathrm{asl}_2(\K)$, the corresponding Lie superalgebra
is the classical simple Lie antialgebra
$\mathrm{osp}(1|2)$.
We will give the explicit construction of this
Lie superalgebra in the next section
and use it as the main tool for our study.

%%%%%%%%%%%%%%%%%%%%%%
\section{Representations of $\mathrm{asl}_2(\K)$ and
the ghost Casimir of $\mathrm{osp}(1|2)$}
%%%%%%%%%%%%%%%%%%%%%%

In this section we collect the general information about the
representations of $\mathrm{asl}_2(\K)$.
We also introduce the action of $\mathrm{osp}(1|2)$ and
prove Theorem \ref{NoOne}.

%%%%%%%%%%%%%%%%%%%%%%
\subsection{Generators of $\mathrm{asl}_2(\K)$ and
the $\Z_2$-grading.}
%%%%%%%%%%%%%%%%%%%%%%

Consider an $\mathrm{asl}_2(\K)$-representation 
$V=V_0\oplus V_1$ with
$\chi:\mathrm{asl}_2(\K) \to\End(V)$.
The homomorphism condition (\ref{eqRep})
can be written explicitly in terms of the basis elements:
$$
\left \lbrace
\begin{array}{rcl}
\chi_a\chi_b-\chi_b\chi_a&=& \half \,\chi_\e\\[5pt]
\chi_a\chi_\e +\chi_\e \chi_a&=& \half \,\chi_a\\ [5pt]
\chi_b\chi_\e +\chi_\e \chi_b&= &\half\, \chi_b\\[5pt]
\chi_\e\chi_\e& =& \half \chi_\e.
\end{array}
\right.
$$
Let us simplify the notations by fixing the following elements of $\End(V)$:
\begin{eqnarray*}
A=2\,\chi_a,
\qquad 
B=2\,\chi_b, 
\qquad
\E=2\,\chi_\e.
\end{eqnarray*}
The above relations read:
\begin{equation} 
\label{systemrep2}
\left \lbrace
\begin{array}{rcl}
AB -BA&=&\E\\ [5pt]
A\E + \E A&=&A\\ [5pt]
B\E + \E B&=&B\\ [5pt]
 \E^2&=&\E.
\end{array}
\right.
\end{equation}

The element $\E$ is a projector in $V$.
This leads to a decomposition of $V$ into eigenspaces
$V=V^{(0)}\oplus V^{(1)}$ defined by
$$
V^{(\lambda)}=
\{v \in V \; | \; \E v =
\lambda\, v \, \},
\qquad
\lambda=0, 1.
$$
This decomposition is not necessarily the same as the initial one, 
$V=V_0\oplus V_1$. 
Since $V_i$, where $i=0,1$, is stable under the action of 
$\E$, this gives a refinement:
$$
V=V_0^{(0)}\oplus V_0^{(1)}\oplus V_1^{(0)}\oplus V_1^{(1)}
$$
where 
$$
V_i^{(\lambda)}=
\{v \in V_i \; | \; \E v =\lambda\, v \, \}, 
\qquad
\lambda=0, 1,
\quad
i=0, 1.
$$

\begin{prop}
\label{Zeddeux}
Any $\mathrm{asl}_2(\K)$-representation $(V=V_0\oplus V_1,\chi)$ is equivalent to a representation $(V'=V'_{0}\oplus V'_{1}, \chi')$ such that 
\begin{equation*}
\E\left|_{V'_0}\right.=0,
\qquad
\E\left|_{V'_1}\right.=\id.
\end{equation*}
\end{prop}

\begin{proof}
Using the relations (\ref{systemrep2}) it is easy to see
that $A$ and $B$ send the spaces $V_i^{(\lambda)}$ 
into $V_{1-i}^{(1-\lambda)}$, where $\lambda=0, 1,\; i=0, 1$. 
Thus, changing the $\Z_2$-grading of $V$ to  $V'=V'_{0}\oplus V'_{1}$ where
\begin{eqnarray*}
V'_0&=&V_0^{(0)}\oplus V_1^{(0)},\\
V'_1&=&V_0^{(1)}\oplus V_1^{(1)},
\end{eqnarray*}
does not change the parity of the operators $A$, $B$ and $\E$ viewed 
as elements of the $Z_2$-graded space $\End(V')$. 
In other words, the map 
$\chi':\ga\to\End(V')$ defined by $\chi'_x=\chi_x$ for all $x\in \ga$, 
is still an even map satisfying the condition \eqref{eqRep}. 
By consequent, $(V',\chi')$ is also a representation. 
It is then clear that  $(V',\chi')$ is equivalent (in the sense of Definition \ref{defRep} (d))
to $(V,\chi)$.
\end{proof}

As a consequence of the above proposition we will always assume in the sequel that the $\Z_2$-grading of the representations $V$ is given by the eigenspaces of the action $\E$.

%%%%%%%%%%%%%%%%
\subsection{Action of  $\mathrm{osp}(1|2)$}
%%%%%%%%%%%%%%%%

This section provides a special case of the general construction of \cite{Ovs},
 For the sake of completeness, we give here the details of the computations.

Given an $\mathrm{asl}_2(\K)$-representation  $V$,
we define the operators $E$, $F$ and $H$ by
\begin{equation}\label{EFH}
E=A^2, 
\quad 
F=-B^2, 
\quad H=-(AB+BA).
\end{equation}
These three elements define a structure of
 $\mathrm{sl}_2(\K)$-module
on $V$.

\begin{lem}
One has:
\begin{equation}
\begin{array}{rcl}
\label{EFH}
[H,E]&=&2E\\[5pt]
[H,F]&=&-2F\\[5pt]
[E,F]&=&H
\end{array}
\end{equation}
\end{lem}

\begin{proof} 
These relations follow from relations (\ref{systemrep2}).
Indeed,
\begin{eqnarray*}
[H,E]&=&-(AB+BA)A^2+A^2(AB+BA)\\
&=&-ABA^2-BA^3+A^3B+A^2BA\\
&=&-2ABA^2+(ABA^2-BA^3)+(A^3B-A^2BA)+2A^2BA\\
&=&-2ABA^2+(AB-BA)A^2+A^2(AB-BA)+2A^2BA\\
&=&2A(AB-BA)A+(AB-BA)A^2+A^2(AB-BA)\\
&=&2A\E A+\E A^2+A^2\E \\
&=&( A\E A +\E A^2)+ (A\E A + A^2\E)\\
&=&(A\E + \E A)A+ A(\E A+ A \E)\\
&=& A^2+A^2\\
&=&2E.
\end{eqnarray*}
In the same way we obtain $[H,F]=-2F$.
Finally,
\begin{eqnarray*}
[E,F]&=&-A^2B^2+B^2A^2\\
&=& -A^2B^2+ABAB-ABAB+BA^2B-BA^2B+BABA-BABA +B^2A^2\\
&=& -A(AB-BA)B-(AB-BA)AB-BA(AB-BA)-B(AB-BA)A\\
&=&-A\E B-\E AB-BA\E-B\E A\\
&=& -(A\E+ \E A)B-B(A\E+\E A)\\
&=&-AB-BA\\
&=&H
\end{eqnarray*}
\end{proof}

With similar computations one can establish the following additional relations:
\begin{equation}
\label{commutations}
\begin{array}{llll}
[H, A]&=A \qquad\qquad &[E,A]=0 \qquad\qquad &[F,A]=B,\\[5pt]
[H, B]&=-B            &[E,B]=A              &[F,B]=0,\\[5pt]
[H, \E]&=0 &[E,\E]=0 &[F,\E]=0,
\end{array}
\end{equation}
that can be summarized as follows.

\begin{prop}
Every representation of $\mathrm{asl}_2(\K)$ has a structure
of a module over the Lie superalgebra $\mathrm{osp}(1|2)$.
 The even part,
$\mathrm{osp}(2|1)_0$, is spanned by $E,F,H$, 
while the odd part, $\mathrm{osp}(2|1)_1$,
is spanned by $A,B$.
\end{prop}

The relation (\ref{EFH}) and (\ref{commutations}) play an essential role
in all our computations.

%%%%%%%%%%%%%%%%%%%%%%%%%%%%%%%%
\subsection{The ghost Casimir element.}\label{PfNoOne}
%%%%%%%%%%%%%%%%%%%%%%%%%%%%%%%%
The notion of \textit{twisted adjoint action} was introduced for a certain
class of Lie superalgebras in \cite{Fra}.
We recall here the definition in the $\mathrm{osp}(1|2)$-case.

Let $X$ be an element of $\mathrm{osp}(1|2)$ and $Y$ be an element of the
universal enveloping algebra $\U(\mathrm{osp}(1|2))$.
Define $\widetilde{\ad}:\mathrm{osp}(1|2)\to\End(\U(\mathrm{osp}(1|2)))$ by
\begin{equation}
\label{Twist}
\widetilde{\ad}_XY:=
XY-(-1)^{p(X)(p(Y)+1)}\,YX.
\end{equation}
In other words, 
$$
\widetilde{\ad}_X=
\left\{
\begin{array}{rl}
\ad_X & \hbox{if $X$ is even}\\[4pt]
-\ad_X & \hbox{if $X$ is odd}.
\end{array}
\right.
$$
Remarkably enough, $\widetilde{\ad}$ defines an $\mathrm{osp}(1|2)$-action
on $\U(\mathrm{osp}(1|2))$.

The \textit{ghost Casimir} elements
are the invariants of the twisted adjoint action,
see \cite{Fra}, and also \cite{Gor}.
In the case of $\mathrm{osp}(1|2)$, the ghost Casimir element
is particularly simple:
\begin{equation}
\label{GCEl}
\textstyle
\cc=AB-BA-\frac{1}{2}\,\id,
\end{equation}
The ghost Casimir satisfies $\widetilde{\ad}_X\cc=0$ for
all $X\in\mathrm{osp}(1|2)$.

The relation between the above twisted adjoint action and our situation is the following.
Consider a representation $\chi:\mathrm{asl}_2(\K)\to\End(V)$. Denote by $U$ the subalgebra of $\End(V)$ generated by the image of $\mathrm{asl}_2(\K)$ under $\chi$. The algebra $U$ can be viewed as a quotient of $\U(\mathrm{osp}(1|2))$. 
A simple comparison of (\ref{Twist}) and (\ref{AntiCoCo}) shows
that, if $x$ is an odd element of $\mathrm{asl}_2(\K)$, then
$$
]\chi_x,Y[=\widetilde{\ad}_{\chi_x}(Y),
$$
for all $Y$ in $U$.
The operator $\E$ and the ghost Casimir $\cc$ are obviously
related by
$$
\textstyle
\cc=\E-\frac{1}{2}\,\id.
$$
It follows that the second and third relations in (\ref{systemrep2})
are equivalent to $\widetilde{\ad}_A\cc=0$ and $\widetilde{\ad}_B\cc=0$,
respectively, while the relation $\E^2=\E$ reads $\cc^2=\frac{1}{4}\,\id$.

This completes the proof of Theorem \ref{NoOne}.

%%%%%%%%%%%%%%%%%%%%%%%%%%%%%%%%
\subsection{Usual Casimir elements.}
%%%%%%%%%%%%%%%%%%%%%%%%%%%%%%%%
The operator $\E$ is also related to the usual Casimir elements $C$, resp. $C_0$  of $\mathrm{osp}(1|2)$, resp. $\mathrm{sl}_2(\K)$. Recall
\begin{eqnarray*}
C&=&\textstyle E{}F+F{}E+
\frac{1}{2}
\left(
H^2+A{}B-B{}A
\right),\\
C_0&=&\textstyle E{}F+F{}E+
\frac{1}{2}
H^2.
\end{eqnarray*}
We easily see
$$
\E=2\left(
C-C_0
\right).
$$
This implies that, if $V$ is an irreducible representation of $\mathrm{asl}_2(\K)$,
then $\E|_{V_0}$ and $\E|_{V_1}$ are proportional to $\id$.

Moreover, straightforward computation in  $\U(\mathrm{osp}(1|2))$ gives the following relation
$$
4(C-C_0)^2=4C-2C_0.
$$
It follows the condition $\E^2=\E$ implies that $C$ acts trivially.

%%%%%%%%%%%%%%%%%%%%%%%%%%%%%%%%
%%%%%%%%%%%%%%%%%%%%%%%%%%%%%%%%
\section{Weighted representations of $\mathrm{asl}_2(\K)$.}
%%%%%%%%%%%%%%%%%%%%%%%%%%%%%%%%
%%%%%%%%%%%%%%%%%%%%%%%%%%%%%%%%

In this section, we introduce the notion of weighted representation
of the Lie antialgebra $\mathrm{asl}_2(\K)$.
This class of representation is characterized by the property that
the action of the Cartan element $H$ of $\mathrm{osp}(1|2)$ has at least one eigenvector.
We do not require \textit{a priori}  the  eigenspaces to be finite dimensional.

%%%%%%%%%%%%%%%%%%%%%%%%%%%%%%%%
\subsection{The definition.}
%%%%%%%%%%%%%%%%%%%%%%%%%%%%%%%%

 Let $V$ be a representation of $\mathrm{asl}_2(\K)$.
We introduce the subspaces
$$ 
V_\ell=\{ v\in V\,|\, Hv=\ell\,v\}, \; \ell \in \K.
$$ 
Whenever $V_\ell \not =\{0\}$, we call this subspace a 
\textit{weight space} of $V$ with weight $\ell$. 
We denote by $\Pi_H(V)$ the set of weights
of representation $V$.

\begin{lem}
With the above notations:
\begin{enumerate}
\item[(i)] The element $A$ (resp. $B$) maps $V_\ell$ into $V_{\ell +1}$ (resp. $V_{\ell-1}$).

\item[(ii)] The sum $\sum _{\ell \in \Pi_H(V)} V_\ell$ is direct in $V$.

\item[(iii)]
The space 
$$
Wt(V):=\bigoplus_{\ell \in \Pi_H(V)} V_\ell
$$
is a subrepresentation of $V$.
\end{enumerate}
\end{lem}

\begin{proof}
Let $v$ be a vector in $V_\ell$.
Using the relations (\ref{commutations}) we obtain:
$$
HAv=[H, A]v+AHv=Av+\ell Av=(\ell+1)v
$$
and
$$
HBv=[H, B]v+BHv=-Bv+\ell Bv=(\ell-1)v.
$$
Part (i) then follows.

Part (ii) is clear since the weight spaces are eigenspaces for $H$. 

It follows from (i) that $Wt(V)$ is stable with respect to the action of
$A$ and $B$ and, therefore,
it is also stable under $\E=AB-BA$.
Hence (iii).
\end{proof}

\begin{cor}
If $V$ is an irreducible representation then either
$$ Wt(V)=\{0\} \quad \text{ or } \quad Wt(V)=V.$$
\end{cor}

\begin{defn}
We call \textit{weighted representation} any representation $V$ of $\mathrm{asl}_2(\K)$ such that $Wt(V)\not=\{0\}$.
\end{defn}

%%%%%%%%%%%%%%%%%%%
\subsection{The family of weighted representations $V(\ell)$.}\label{Vell}
%%%%%%%%%%%%%%%%%%%

For every $\ell \in\K$, we construct an irreducible weighted
representation of $\mathrm{asl}_2(\K)$ that we denote $V(\ell)$.
This representation contains an odd vector $e_1$ such that $He_1=\ell e_1$
and, by irreducibility, every element of $V(\ell)$ is a result of the
(iterated) $\mathrm{asl}_2(\K)$-action on $e_1$.

\medskip

(a)
\textbf{The case where $\ell$ is not an odd integer.}
We start the construction with the generic weight $\ell$.

Consider a family of linearly independent vectors  $\{e_k \}_{k\in \Z}$. 
We set  $\displaystyle V(\ell)=\bigoplus_{k\in \Z} \K e_k$ and we
define the operators $A$ and $B$ on $V(\ell)$ by
\begin{eqnarray*}
Ae_k&=&e_{k+1}, \; \forall k\in \Z \\[5pt]
Be_k&=& \big( (1-\ell)/2-[k/2]\big)e_{k-1} , \; \forall k\in \Z, 
\end{eqnarray*}
The operator $\E$ is determined by $\E=AB-BA$.
Introduce the following $\Z_2$-grading on $V(\ell)$:
\begin{equation}
\label{GradEq}
\begin{array}{rcl}
V(\ell)_0&=&
\displaystyle 
\bigoplus_{k \text{ even}}\K e_k,\\[14pt]
V(\ell)_1&=&
\displaystyle
\bigoplus_{k \text{ odd}}\K e_k.
\end{array}
\end{equation}
It is easy to see the operators $A$ and $B$ are odd operators 
with respect to this grading whereas $\E$ is even.

\begin{prop}
\label{DefRepLem}
The space $V(\ell)$ together with the operators $A, B, \E$ 
is an $\mathrm{asl}_2(\K)$-representation.
\end{prop}

\begin{proof}
By simple straightforward computations  we obtain:
$$
A\E  +\E A=A,
\qquad
B\E +\E B=B.
$$
Moreover,
on the basis elements $e_k$ of $V(\ell) $
\begin{eqnarray*}
\E e_k&=&
\begin{cases} e_k, &\hbox{if $k$ is odd}\label{r'2}\\[5pt]
0, &\hbox{if $k$ is even} \end{cases}\\
%He_k&=&(l+k)e_k \label{r3}
\end{eqnarray*}
so that $\E^2=\E$.
\end{proof}

It is easy to see the basis elements $e_k$'s are weight vectors.
Indeed, one checks
\begin{equation}
\label{WeightVect}
He_k=(\ell+k-1)\,e_k, \qquad \forall\; k\in\Z.
\end{equation}
In particular, the element $e_1$ is a weight vector of weight 
$\ell$ and generates the representation $V(\ell)$.

The actions on the basis elements of $\mathrm{osp}(1|2)$ can be pictured as follows:
\medskip

 \SelectTips{eu}{12}%
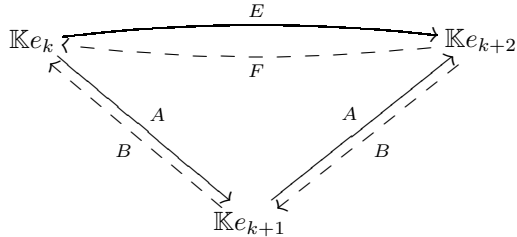
\xymatrix{&&& \K e_k \ar@/^/[rrrr]^{E} \ar@
<2pt>[rrdd]^A&&&& \K e_{k+2}\ar@/^/@{-->}[llll]^{F}\ar@<2pt>@{-->}[lldd]^B\\
\\
&&&	&& \K e_{k+1}\ar@<2pt>@{-->}[lluu]^B \ar@<2pt>[rruu]^A}

\medskip

\noindent
The entire space $V(\ell)$ can be pictured as a infinite chain of 
the above diagrams.

\medskip

 
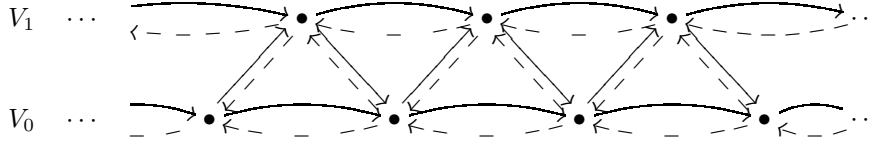
\xymatrix{
  V_1 \quad \cdots\quad \ar@/^/[rr]  && 
  \bullet \ar@/^/[rr] \ar@/^/@{-->}[ll]\ar@<2pt>@{-->}[ld]  
  \ar@<2pt>[rd]&& \bullet \ar@/^/[rr]\ar@/^/@{-->}[ll] 
  \ar@<2pt>@{-->}[ld]  \ar@<2pt>[rd]&& 
  \bullet  \ar@<2pt>@{-->}[ld]   \ar@<2pt>[rd]\ar@/^/@{-->}[ll] \ar@/^/[rr]&& 
  \ar@/^/@{-->}[ll] \cdots\\ 
V_0  \quad \cdots\quad \ar@/^/[r]    & 
\bullet \ar@/^/@{--}[l]  \ar@/^/[rr] \ar@<2pt>[ru] && 
\bullet \ar@<2pt>[ru]\ar@<2pt>@{-->}[lu] \ar@/^/[rr] \ar@/^/@{-->}[ll] &&
 \bullet\ar@<2pt>@{-->}[lu]\ar@<2pt>[ru] \ar@/^/@{-->}[ll] \ar@/^/[rr]
&&\ar@<2pt>@{-->}[lu] \ar@/^/@{-->}[ll] \bullet  \ar@/^/@{-}[r]& \ar@/^/@{-->}[l]\cdots}

\bigskip 

(b)
\textbf{Construction of $V(\ell)$ for $\ell$ a positive odd integer.}
Consider a family of linearly independent vectors  $\{e_k \}_{k\in \Z,k\geq 2-\ell}$.
We set  $\displaystyle V(\ell)=\bigoplus_{k\geq 2-\ell} \K e_k$ and
we define the operators $A$ and $B$ on $V(\ell)$ by  a similar formula:
\begin{eqnarray*}
Ae_k&=&e_{k+1}, \; \forall k\geq 2-\ell\\[5pt]
Be_k&=& \big( (1-\ell)/2-[k/2]\big)e_{k-1} , \; \forall k>2-\ell, \\[5pt]
Be_{2-\ell}&=&0.
\end{eqnarray*}
The operator $\E$ is again determined by $\E=AB-BA$.
The $\Z_2$-grading on $V(\ell)$ is defined by the same
formula (\ref{GradEq}).
The result of Lemma \ref{DefRepLem} holds true.

 The element $e_1$ is an odd weight vector of weight $\ell$ 
 and generates the representation $V(\ell)$. 
 However, the vector $e_{2-\ell}$ is more interesting. 
 
 \begin{defn}
 \label{LowwLem}
 We call a representation $V$ a lowest weight  (resp. highest weight)
 representation if it contains a vector
 $v$, such that $Bv=0$ and the vectors $A^nv$ span $V$
 (resp. $Av=0$ and $B^v$ span $V$);
 the vector $v$ is called a lowest weight (resp. highest weight) vector.
  \end{defn}
  
Clearly, the vector $e_{2-\ell}$ is a lowest weight vector of the representation
$V(\ell)$, if $\ell$ a positive odd integer.
One obtains the following diagram.

\bigskip
 
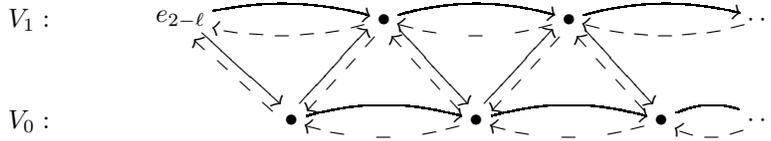
\xymatrix{
  V_1 :\quad  &e_{2-\ell} \ar@/^/[rr]  \ar@<2pt>[rd]&& \bullet \ar@/^/[rr]\ar@/^/@{-->}[ll] \ar@<2pt>@{-->}[ld]  \ar@<2pt>[rd]&& \bullet  \ar@<2pt>@{-->}[ld]   \ar@<2pt>[rd]\ar@/^/@{-->}[ll] \ar@/^/[rr]&& \ar@/^/@{-->}[ll] \cdots\\ 
V_0  :\quad  && \bullet \ar@<2pt>[ru]\ar@<2pt>@{-->}[lu] \ar@/^/[rr] && \bullet\ar@<2pt>@{-->}[lu]\ar@<2pt>[ru] \ar@/^/@{-->}[ll] \ar@/^/[rr]
&&\ar@<2pt>@{-->}[lu] \ar@/^/@{-->}[ll] \bullet  \ar@/^/@{-}[r]& \ar@/^/@{-->}[l]\cdots}
\bigskip

Viewed as a representation of $\mathrm{osp}(1|2)$, $V(\ell)$ is a Verma module.

\medskip

(c)
 \textbf{Construction of $V(\ell)$ for $\ell$ a negative odd integer.}
Consider a family of linearly independent vectors  $\{e_k \}_{k\in \Z,k\leq -\ell}$.
We set  $\displaystyle V(\ell)=\bigoplus_{k\leq -\ell} \K e_k$ and we define the operators $A$ and $B$ on $V(\ell)$ by
\begin{eqnarray*}
Ae_k&=& \big((1+\ell)/2+[k/2]\big)e_{k+1} , \; \forall k<-\ell, \\[5pt]
Ae_{-\ell}&=&0,\\[5pt]
Be_k&=&e_{k-1}, \; \forall k\leq -\ell .\\[5pt]
\end{eqnarray*}
As previously these operators define an $\mathrm{asl}_2(\K)$-representation. 
The vector $e_{-\ell}$ is a highest weight vector of $V(\ell)$.

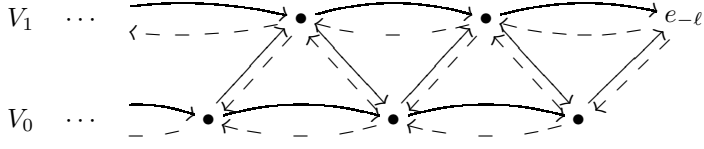
\medskip

 \xymatrix{
  V_1 \quad \cdots\quad \ar@/^/[rr]  && \bullet \ar@/^/[rr] \ar@/^/@{-->}[ll]\ar@<2pt>@{-->}[ld]  \ar@<2pt>[rd]&& \bullet \ar@/^/[rr]\ar@/^/@{-->}[ll] \ar@<2pt>@{-->}[ld]  \ar@<2pt>[rd]&& e_{-\ell}\ar@<2pt>@{-->}[ld]\ar@/^/@{-->}[ll] \\ 
V_0  \quad \cdots\quad \ar@/^/[r]    & \bullet \ar@/^/@{--}[l]  \ar@/^/[rr] \ar@<2pt>[ru] && \bullet \ar@<2pt>[ru]\ar@<2pt>@{-->}[lu] \ar@/^/[rr] \ar@/^/@{-->}[ll] && \bullet\ar@<2pt>@{-->}[lu]\ar@<2pt>[ru] \ar@/^/@{-->}[ll] }

\bigskip

%%%%%%%%%%%%%%%%%%%%%%
\subsection{Geometric realization.}
%%%%%%%%%%%%%%%%%%%%%%
It was shown in \cite{Ovs} that $\mathrm{asl}_2(\K)$ has a representation
in terms of vector fields on the ${1|1}$-dimensional space.
More precisely, consider $\mathscr{F}=C^{\infty}_\K(\R)$ the set of $\K$-valued $C^{\infty}$-functions of one real variable $x$. Introduce $\mathscr{A}=\mathscr{F}[\xi]/ (\xi^2)$ with the $\Z_2$-grading $ \mathscr{A}_0=\mathscr{F}$, $\mathscr{A}_1=\mathscr{F}\xi$. Define the vector field
$$
\Dc=\frac{\partial}{\partial\xi}+\xi\,\frac{\partial}{\partial{}x}.
$$
It is very easy to check that the following vector fields:
\begin{equation}
\label{GeomRep}
A=\Dc,
\qquad
B=x\,\Dc,
\qquad
\E=\xi\,\Dc
\end{equation}
satisfy the relations \eqref{systemrep2}.
Therefore this defines an $\mathrm{asl}_2(\K)$-action on $\mathscr{A}$.

Let us choose the following function:
$$
e_1=x^\l\,\xi,
$$
where $\l\in\K$.
It turns out this function generates a weighted representation of
$\mathrm{asl}_2(\K)$, isomorphic to $V(\ell)$, with the weight
$\ell=-2\l-1$. Note that the case $\lambda$ is an integer gives the highest or lowest irreducible representation.

\begin{rem}
The vector field $\Dc$ defines the standard contact structure on $\K^{1|1}$.
The vector fields (\ref{GeomRep}) are therefore \textit{tangent} to the contact structure.
These vector fields do not span a Lie (super)algebra with respect to the usual
Lie bracket. 
Remar\-kable enough, the corresponding generators of
the $\mathrm{osp}(1|2)$-action:
$$
E=\frac{\partial}{\partial{}x},
\qquad
F=-x^2\,\frac{\partial}{\partial{}x}-\xi\,\Dc,
\qquad
H=-2x\,\frac{\partial}{\partial{}x}-\xi\,\Dc
$$ 
are contact vector fields, while the above vector fields
$A$ and $B$ are the only vector fields
that are contact and tangent at the same time.
The Lie superalgebra $\mathrm{osp}(1|2)$ thus preserves the contact structure
in the usual way.
\end{rem}

%%%%%%%%%%%%%%%%%%%%%%%%%
%%%%%%%%%%%%%%%%%%%%%%%%%
\section{Classification results}
%%%%%%%%%%%%%%%%%%%%%%%%%
%%%%%%%%%%%%%%%%%%%%%%%%%

In this section we prove Theorem \ref{InfThm} which is our main classification result.

%%%%%%%%%%%%%%%%%%%%%%
\subsection{Classification of weighted representations.}
%%%%%%%%%%%%%%%%%%%%%%
The following statement shows that the representations $V(\ell)$ are,
indeed, irreducible; it also classifies all the isomorphisms between these
representations.

\begin{prop}
\label{MainThm}
(i)
For every $\ell\in\K$, the representation $V(\ell)$ is an irreducible
infinite-dimensional representation.

(ii)
Two weighted representations $V(\ell)$ and $V(\ell')$,
where $\ell$  and $\ell'$ are not odd integers, are isomorphic if and only if 
$\ell'-\ell=2\,m$ for some $m\in\Z$.

(iii)
Two weighted representations $V(\ell)$ and $V(\ell')$,
where $\ell$  and $\ell'$ are odd integers, 
are isomorphic if and only if $\ell$ and $\ell'$ have same sign.
\end{prop}

\begin{proof}
(i)
Irreducibility of $V(\ell)$:

Suppose there exists a subrepresentation $V'$ of $V(\ell)$. Consider a nonzero vector $v\in V'$ and write 
$$
v=\displaystyle \sum_{1\leq i \leq N} \alpha_i e_{k_i}
$$
with $\alpha_i\not =0$, for all  $1\leq i \leq N$.
Using (\ref{WeightVect}), we obtain,
$$
H^pv=\displaystyle \sum_{1\leq i \leq N} \alpha_i (\ell+k_i-1)^pe_{k_i}, 
$$
where $0\leq p\leq N-1$.
We may assume the coefficients $(\ell+k_i-1)$ are nonzero.
In other words, we assume the basis elements occurring in the decomposition of 
$v$ are not of weight zero. 
If this is not the case, we change $v$ to $A^kv$ or $B^kv$ for sufficiently large $k$.

The above equations form a linear system of type Vandermonde 
with distinct nonzero coefficients, and so it is solvable. 
We can express the $e_{k_i}$'s  as linear combinations of 
$H^pv$'s. 
It follows that the vectors $e_{k_i}$ are in $V'$. 
Applying $A$ and $B$ to $e_{k_i}$, we produce all the vectors $e_k$, $k \in \Z$. 
Hence, the vectors $e_k$, of the basis are in $V'$ for all $k \in \Z$.
This implies $V'=V(\ell)$. 
Therefore  there is no proper subrepresentation of $V(\ell)$.

(ii) Let $\ell$ and $\ell'$ be two scalars in $\K$ that are not odd integers.

Denote by $\{e_k\}_k$ the standard basis of $V(\ell)$ and $\{e'_{k}\}_k$ the standard basis of $V(\ell')$. Suppose there exists an isomorphism of representation $\Phi: V(\ell') \mapsto V(\ell)$.

The vector $\Phi(e'_1)$ is a vector of weight $\ell'$ in $V(\ell)$. 
The weights in $V(\ell)$ 
are of the form $\ell+k$ for some $k\in \Z$ and the corresponding weight space is $\K e_{k+1}$. 
We thus have $\Phi(e'_1)=\a\,e_{k+1}$ for some $\a \in \K \setminus \{0\}$,  
$k\in \Z$, and therefore $\ell'=\ell+k$.

Moreover, $\E\,\Phi(e'_1)=\Phi(e'_1)$, so $e_{k+1}$ has to be an odd vector, $i.e.$, $k$ has to be even. 

We proved that a necessary condition to have $V(\ell)$ isomorphic to $V(\ell')$ is
$$
\ell'= \ell+2m,
$$
where $m\in \Z$.

Conversely, suppose $\ell'= \ell+2m,$ for some $m\in \Z$. 
It is easy to check that the linear map  $\Phi: V(\ell') \mapsto V(\ell)$ defined by 
$$
\Phi(e'_k)=e_{k+2m},
$$
for all $k\in \Z$,
is an isomorphism of representation.

(iii) 
Let $\ell$ and $\ell'$ be two odd integers.
If $\ell$ and $\ell'$ have opposite sign then $V(\ell)$ and $V(\ell')$ 
cannot be isomorphic, since on one of the space $A$ acts injectively and on the other space $A$ has a non-trivial kernel.

Conversely, if $\ell$ and $\ell'$ have same sign, let us construct an explicit isomorphism between 
$V(\ell)$ and $V(\ell')$. 
One has:  $\ell'= \ell-2m,$ for some $m\in \N$ and we define  $\Phi: V(\ell') \mapsto V(\ell)$ by 
$$
\Phi(e'_k)=e_{k+2m}, \quad \forall k\geq 2-\ell'
$$
in the case where $\ell$ and $\ell'$ are positive, or by
$$
\Phi(e'_k)=e_{k+2m}, \quad \forall k\leq -\ell'
$$
in the case where $\ell$ and $\ell'$ are negative.
\end{proof}

%%%%%%%%%%%%%%%%%%%%%%
\subsection{Structure of weighted representations}
%%%%%%%%%%%%%%%%%%%%%%

In this section we establish two lemmas crucial 
for the proofs of Theorem \ref{FinThm} and \ref{InfThm}. 

Let us consider an irreducible representation $V$ of $\mathrm{asl}_2(\K)$.
Recall that the $\Z_2$-grading $V=V_0\oplus V_1$ corresponds to
the decomposition with respect to the eigenvalues of $\E$,
see Lemma \ref{Zeddeux}.

\begin{lem} \label{lem1}
If there exists a nonzero element $v \in V$ such that $Hv=\ell v$, $\ell \in \K$,
then there exist nonzero elements $v_i\in V_i$  and $\ell_i \in \K$ , $i=0,1$, 
such that  $Hv_i=\ell_iv_i$.
\end{lem}

\begin{proof}
In the case where $v$ is a pure homogeneous element, $i.e.$, $v\in V_i$, we can choose
for $v_{1-i}$ the vector $Av$ or $Bv$.  Indeed, one has
 \begin{equation}\label{HAHB}
\begin{array}{lcl}
HAv&=&[H,A]v+AHv=Av+\ell Av=(\ell+1)v,\\[5pt]
HBv&=&[H,B]v+BHv=-Bv+\ell Av=(\ell-1)v.
\end{array}
\end{equation}
If the two vectors $Av$ and $Bv$ are zero then $V=\K v$
is the trivial representation and, by consequent, the statement of the theorem is true. 
Otherwise, if $Av$ and $Bv$ are not zero,
one obtains weight vectors of weight $\ell_{1-i}=\ell{}\pm 1$. 

If $v$ is not a homogeneous element, we can write $v=v_0+v_1$ with
$v_0\not =0 \in V_0$ and $v_1\not= 0 \in V_1$.
One then has: 
$$Hv=Hv_0+Hv_1 \text{ and }Hv=\ell v=\ell v_0+\ell v_1.$$
Furthermore, $Hv_0$ is an element of $V_0$ and $Hv_1$ is an element of $V_1$,
since the operator $H=-(A\,B+B\,A)$ is even.
Therefore, by uniqueness of the writing in $V_0\oplus V_1$, one has:
$$Hv_0=\ell v_0 \text{ and }Hv_1=\ell{}v_1.$$
\end{proof}

\begin{lem}
\label{lambdamu}
If $v \in V_i$ is such that $Hv=\ell{}v$ then :
\item[(i)] for all $k\geq 1$,
$$ AB^kv=\lambda_kB^{k-1}v$$
with $\lambda_k = \displaystyle \Big[ \dfrac{k-i}{2}\Big] +\dfrac{i-\ell{}}{2}$,
\noindent
where $ \Big[ \dfrac{k-i}{2}\Big] $ is the integral part of $(k-i)/2$;

\medskip

\item[(ii)] for all $k\geq 1$,
$$ BA^kv=\mu_kA^{k-1}v$$
with $\mu_k = -\displaystyle \Big[ \dfrac{k-i}{2}\Big] -\dfrac{i+\ell{}}{2}$.

\end{lem}

\begin{proof}
We first establish the formulas for $k=1$.
On the one hand one has:
$$ ABv=-Hv-BAv=-\ell{}v-BAv$$
and on the other hand,
$$ABv=\E v +BAv= iv+BAv.$$
By adding or subtracting these two identities we deduce:
$$ 2ABv=(i-\ell{})v \text{ and } 2BAv=(-i-\ell{})v.$$
We get $\lambda_1=(i-\ell{})/2$ and $\mu_1=(-i-\ell{})/2$,
so that the formulas are established at the order 1. 

By induction on $k$,
\begin{eqnarray*}
ABB^{k-1}v&=&-HB^{k-1}v-BAB^{k-1}v \\
&=& -(\ell{}-(k-1))B^kv - \lambda_{k-1}B^{k-1}v\\
&=& (-\ell{}+k-1- \lambda_{k-1})B^{k-1}v
\end{eqnarray*}
We deduce the relations:
\begin{eqnarray*}
\lambda_{k}&=&-\ell{}+k-1- \lambda_{k-1}\\
&=&-\ell{}+k-1-(-\ell{}+(k-2)-\lambda_{k-2})\\
&=& 1+ \lambda_{k-2}.
\end{eqnarray*}
Knowing $\lambda_1=\dfrac{i-\ell{}}{2}$ we can
now obtain the explicit expression of $\lambda_k$:
$$\lambda_k = \displaystyle \Big[ \dfrac{k-i}{2}\Big] +\dfrac{i-\ell{}}{2}.$$
Hence part (i).

Part (ii) can be proved in a similar way.
\end{proof}

%%%%%%%%%%%%%%%%%%%%%%%%
\subsection{Proof of Theorem \ref{FinThm}}
%%%%%%%%%%%%%%%%%%%%%%%%

In this section we prove 
that  $\mathrm{asl}_2(\K)$ has no non-trivial representations of finite dimension.

Let $V$ be an irreducible finite dimensional representation of $\mathrm{asl}_2(\K)$. 
Considering the actions of the elements $E=A^2$, $F=-B^2$ and $H$, 
the space $V$ has a structure of $\mathrm{sl}_2$-module. 
Therefore, there exists a weight vector $v$ such that $Hv=\ell\,v$ 
for some $\ell \in \Z$.

By lemma \ref{lem1}, we can assume that $v$ is a homogeneous element, 
namely $v\in V_i$, $i=0,1$. 
Let us consider the family of vectors
$$
\F=\left\{
\cdots, B^kv, \cdots, Bv, v, Av, \cdots, A^kv, \cdots
\right\}.
$$
From formula (\ref{HAHB}) we know that all the nonzero vectors of $\F$ 
are eigenvectors of $H$ with distinct eigenvalues, $\ell \pm k$, $k\in \N$. 
Therefore, all the non-zero vectors of $\F$ are linearly independent. 

Hence, there exists $N\geq 1$ such that
$$ B^{N-1}v\not=0,  \quad B^{k}v= 0,\; \forall k\geq N$$
and  $M\geq 1$ such that
$$A^{M-1}v\not=0, \quad A^{k}v= 0, \;\forall k\geq M.$$
 Using Lemma \ref{lambdamu} we deduce
\begin{eqnarray*}
\lambda_N&=&0,\\
\mu_M&=&0.
\end{eqnarray*}
This leads to
\begin{eqnarray*}
\displaystyle \Big[ \dfrac{N-i}{2}\Big] +\dfrac{i-\ell}{2}&=&0,\\
-\displaystyle \Big[ \dfrac{M-i}{2}\Big] -\dfrac{i+\ell}{2}&=&0.
\end{eqnarray*}
By subtracting these two equations we obtain:
$$ \Big[ \dfrac{N-i}{2}\Big] +\Big[ \dfrac{M-i}{2}\Big] +i=0.$$
But one has: $N,M\geq 1$ and $i=0,1$. 
So, necessarily,
$$N=M=1, \quad i=0.$$
In conclusion, $v$ is an even vector such that $Av=Bv=0$ and $V$ 
is nothing but the trivial representation.

Finally, if $V$ is an arbitrary finite-dimensional representation,
then $V$ is completely reducible.
This immediately follows from the classical theorem in the
$\mathrm{osp}(1|2)$-case.

Theorem \ref{FinThm} is proved.

%%%%%%%%%%%%%%%%%%%%%%
\subsection{Proof of Theorem \ref{InfThm}} \label{PfInfThm}
%%%%%%%%%%%%%%%%%%%%%%

Let us consider an infinite-dimensional irreducible weighted representation $V$. 

We start by studying the cases of highest weight representations and lowest weight representations.

\begin{lem}\label{highrep}
Every irreducible highest weight representation is isomorphic to $V(-1)$.
\end{lem}

\begin{proof} Consider an irreducible representation $V$ containing a weight vector $v$ of weight $\ell$ such that $Av=0$. 
We write $v=v_0+v_1$ with $v_i\in V_i$, $i=0,1$. 
We also have $Av_i=0$ and $Hv_i=\ell v_i$ for $i=0,1$.

We first show that $v_0=0$. 
Consider the action of $\mathrm{sl}_2(\K)$ on $V_0$,
the vector $v_0$ is a highest weight vector for this action. 
Since $Av_0=0$, we get 
$$
Hv_0=ABv_0=-\E v_0=0.
$$
The vector $v_0$ is a highest weight vector of weight 0 for the action of  $\mathrm{sl}_2(\K)$. 
By consequence, $v_0$ is also a lowest weight vector, \textit{i.e.}, $B^2v_0=0$. 
Thus, the space $\mathrm{Span}\,( v_0, Bv_0)$ is stable under the action of $A$ and $B$. 
Since $V$ is an infinite-dimensional irreducible representation one necessarily has $v_0=0$.

We can assume now that $v$ belongs to $V_1$. 
Let us use Lemma \ref{lambdamu}, part (ii). 
From the relation $BAv=0$, we deduce $\mu_1=0$ and thus $\ell=-1$. 
This implies that all the constants $\lambda_k$ from Lemma \ref{lambdamu}, part (i)  are non-zero. 
By induction we deduce, using Lemma \ref{lambdamu}, part (i), that all the vectors $B^{k}v$, 
$k\in \N$ are non-zero. 
Moreover, these vectors are linearly independent since they are eigenvectors for $H$ 
associated to distinct eigenvalues. 
By setting
$$e_k=B^{1-k}v, \qquad
k\in \Z,\; k\leq 1,$$
we obtain $V(-1)$ as a subrepresentation of $V$.  
We then deduce from the irreducibility assumption that $V\simeq V(-1)$.
\end{proof}

\begin{lem}\label{lowrep}
Every irreducible lowest weight representation is isomorphic to $V(1)$.
\end{lem}

\begin{proof}
Similar to the proof of Lemma \ref{highrep}.
\end{proof}

Now we are ready to prove Theorem \ref{InfThm}.

Fix a weight vector $v\in V_1$ (such a vector exists by Lemma \ref{lem1}) of some weight $\ell\in \K$. 
Consider the family 
$$\F:=\left\{A^kv, \; B^kv, \;k\in \N \right\}.$$

(a)
Suppose that $\ell$ is not an odd integer. 
It is easy to see that the constants $\lambda_k$ and $\mu_k$, $k\in \N$, 
from Lemma \ref{lambdamu} never vanish. 
Indeed,
\begin{equation}\label{oddell}
\begin{array}{lcl}
\lambda_k=0&\Rightarrow& \ell =2\Big[\dfrac{k-1}{2}\Big]+1,\\[9pt]
\mu_k=0&\Rightarrow& \ell=-2\Big[\dfrac{k-1}{2}\Big]-1.
\end{array}
\end{equation}
By induction, we deduce that the elements in $\F$ are different from zero.  
Moreover, the elements of $\F$ are eigenvectors for the operator $H$ 
with distinct eigenvalue $\ell\pm k$, where $k\in \N$, 
so that, they are linearly independent. 

By setting
$$e_k=\begin{cases} A^{k-1}v\; , & k\geq 1 \\[4pt]
B^{1-k}v\;, & k\leq 0
\end{cases}$$
we see that $V(\ell)$ as a subrepresentation of $V$ . 
Again, by irreducibility assumption, we deduce
$V\simeq V(\ell)$.

Finally, using Proposition \ref{MainThm}, part (ii), one has:
$$V\simeq V(\ell'),$$
where $\ell'$ is the unique element of 
$$\mathcal{P}^+=
[-1,1]\cup \left\{
\ell \in \C\; |\; -1\leq \mathrm{Re}(\ell) <1
\right\},$$ 
such that $\ell'=\ell+2m$ for some $m\in\Z$.

(b)
Suppose that $\ell$ is a positive odd integer. 
From the first statement of \eqref{oddell}, 
we deduce the existence of an integer $N\geq 1$,
 such that $\lambda_N=0$ and $\lambda_k\not=0$ for all $k<N$. 
 Hence,
$$ B^{k}v\not=0,\,\;\forall \; k<N,   \quad AB^{N}v= 0.$$
If $B^Nv\not =0$ then this vector is a highest weight vector. 
By Lemma \ref{highrep}, we obtain $V\simeq V(-1)$. 
But, in the highest weight representation $V(-1)$, the set of weights 
is the set of negative integers. 
We obtain a contradiction since $v$ has a positive weight.

It follows that $B^Nv=0$ and
this implies that $B^{N-1}$ is a lowest weight vector. 
Using Lemma~ \ref{lowrep}, we conclude
$$V\simeq V(1).$$

(c)
Suppose finally that $\ell$ is a negative odd integer. 
Then similar arguments show:
$$V\simeq V(-1).$$

Theorem \ref{InfThm} is proved.

\begin{rem}
We proved that any irreducible weighted representation is a Harish-Chandra irreducible representation (\textit{i.e.} the weight spaces are all finite dimensional). A classification of  Harish-Chandra irreducible representations of $\mathrm{osp}(1|2)$ over the field of complex numbers is given in \cite{BP}.  The correspondence between the representations $V(\ell)$ and the representations given in Theorem 5.13 of \cite{BP} is the following:

(a) If $\ell \in \mathcal{P}^+$ is not an odd integer then
$$
V(\ell)\simeq \mathscr{D}(l,\lambda_0),
$$
for the choice $l=0$ and $\lambda_0=\ell / 2$.

(b) The lowest weight representation is 
$$
V(1)\simeq [\lambda_0] \downarrow,
$$
for the unique choice $\lambda_0=-1/2$.

(c) The highest weight representation is 
$$
V(-1)\simeq [\lambda_0] \uparrow,
$$
for the unique choice $\lambda_0=1/2$.
\end{rem}

%%%%%%%%%%%%%%%%%%%%%%
%%%%%%%%%%%%%%%%%%%%%%
\section*{Appendix: tensor product of two representations}
%%%%%%%%%%%%%%%%%%%%%%
%%%%%%%%%%%%%%%%%%%%%%

Given two representations, $V$ and $W$ of
$\mathrm{asl}_2(\K)$, to what extent their tensor product,
$V\otimes{}W$ is again an $\mathrm{asl}_2$-representation?
This question is non-trivial since $\mathrm{asl}_2(\K)$ is not a
Lie algebra.
We will show that $\mathrm{asl}_2(\K)$ does not act on $V\otimes{}W$.
An attempt to define such an action leads to a deformation
of the $\mathrm{asl}_2$-relations by the Casimir
element of $\U(\mathrm{osp}(1|2))$.
The algebraic meaning of this deformation is not yet clear.

The operators $A$ and $B$ have canonical lifts to
$V\otimes{}W$ according to the Leibniz rule:
$$
\widetilde{A}=A\otimes\id+\id\otimes{}A,
\qquad
\widetilde{B}=B\otimes\id+\id\otimes{}B,
$$
since they belong to the $\mathrm{osp}(1|2)$-action.
It is then natural to define the lift of operator $\E$ by
$
\widetilde{\E}:=
\widetilde{A}\widetilde{B}-\widetilde{B}\widetilde{A}.
$
One immediately obtains the explicit formula
$$
\widetilde{\E}=\E\otimes\id+\id\otimes\E+
2\left(
A\otimes{}B-B\otimes{}A
\right).
$$

The following statement is straightforward.

\begin{prop}
The operators  $\widetilde{A},\widetilde{B}$ and $\widetilde{\E}$
satisfy the following relations:
\begin{equation}
\label{SqrEps}
\begin{array}{rcl}
\widetilde{A}\widetilde{B}-\widetilde{B}\widetilde{A}&=&
\widetilde{\E}\\[5pt]
\widetilde{A}\widetilde{\E}+\widetilde{\E}\widetilde{A}&=&
\widetilde{A}\\[5pt]
\widetilde{B}\widetilde{\E}+\widetilde{\E}\widetilde{B}&=&
\widetilde{B}\\[5pt]
\widetilde{\E}^2&=&\widetilde{\E}+4\,\bar{C},
\end{array}
\end{equation}
where
$$
\textstyle
\bar{C}=
E\otimes{}F+F\otimes{}E+
\frac{1}{2}\left(
H\otimes{}H+A\otimes{}B-B\otimes{}A
\right).
$$
\end{prop}

This means that two of the relations (\ref{systemrep2})
are satisfied, but not the last $\mathrm{asl}_2$-relation
$\E^2=\E$.

Let us recall that the element $C\in{}\U(\mathrm{osp}(1|2))$
given by
$$
\textstyle
C=E{}F+F{}E+
\frac{1}{2}
\left(
H^2+A{}B-B{}A
\right)
$$
is nothing but the classical Casimir element.
The operator $\bar{C}$ is the diagonal part of the
standard lift of $C$ to $V\otimes{}W$.
In particular, the operator $\bar{C}$ commutes
with the action of $\mathrm{asl}_2(\K)$ and $\mathrm{osp}(1|2)$:
$$
[\bar{C},\widetilde{A}]=
[\bar{C},\widetilde{B}]=
[\bar{C},\widetilde{\E}]=0.
$$
This is how the Casimir operator of $\mathrm{osp}(1|2)$
appears in the context of representations of $\mathrm{asl}_2(\K)$.

The relations (\ref{SqrEps}) look like a ``deformation'' of the
$\mathrm{asl}_2$-relations (\ref{systemrep2})
with one parameter that commutes with all the generators.
It would be interesting to find a precise algebraic sense of
this deformation.

%%%%%%%%%%%%%%%%%

\end{document}